\documentclass[10pt,a4paper,leqno]{amsart}

\usepackage{mathrsfs}  
\usepackage[dvips]{graphicx}
\usepackage[cmtip,arrow]{xy}
\usepackage{pb-diagram,pb-xy}

\newcommand{\N}{\mathbb{N}}

\newcommand{\pr}{\textrm{pr}}
\newcommand{\M}{\mathscr{M}}
\newcommand{\id}{\textrm{id}}

\newtheorem{thm}[equation]{Theorem}
\newtheorem{lem}[equation]{Lemma}
\newtheorem{prop}[equation]{Proposition}
\newtheorem{cor}[equation]{Corollary}
\newtheorem{defnA}[equation]{Definition}
\newtheorem{exA}[equation]{Example}

\newenvironment{defn}{\begin{defnA}\rm}{\end{defnA}}

\newtheorem{remA}[equation]{Remark}
\newenvironment{rem}{\begin{remA}\rm}{\end{remA}}

\numberwithin{equation}{section}

\title{Equivariant embedding of metrizable $G$-spaces in linear $G$-spaces}
\author{Aasa Feragen}
\begin{document}

\begin{abstract}
Given a Lie group $G$ we study the class $\M$ of proper metrizable $G$-spaces with metrizable orbit spaces, and show that any $G$-space $X \in \M$ admits a closed $G$-embedding into a convex $G$-subset $C$ of some locally convex linear $G$-space, such that $X$ has some $G$-neighborhood in $C$ which belongs to the class $\M$. As corollaries we see that any $G$-ANE for $\M$ has the $G$-homotopy type of some $G$-CW complex and that any $G$-ANR for $\M$ is a $G$-ANE for $\M$.
\end{abstract}
 
\maketitle

\section{Introduction}

In this paper we study spaces in the class $\M$ of proper metrizable $G$-spaces with metrizable orbit spaces, where $G$ is a Lie group. In the classical theory of retracts the Wojdys{\l}awski embedding theorem (see \cite{hu}) ensures that any metrizable space can be embedded as a closed subset of a convex subspace of some Banach space. In the equivariant case S.~Antonyan has proved that for any topological group $G$, any $G$-space $X$ with a $G$-invariant metric admits a $G$-embedding as a closed $G$-subset of a convex $G$-subspace $C$ of some Banach $G$-space $B$ with a $G$-invariant norm \cite[Proposition~8 and Theorem~1]{an*}. 

However, when we are considering extensions and retractions in the class $\M$ with a noncompact Lie group action, then we need to know something about the metrizability of the orbit space $B/G$, or at least of the orbit space of some $G$-neighborhood of $X$ in $C$.

E.~Elfving \cite{e1} has shown that for a linear Lie group $G$, any Palais proper metrizable $G$-space $X$ which is locally compact, separable and finite-dimensional, and which has finitely many orbit types, admits a closed $G$-embedding in a linear $G$-space such that $G$ acts Palais properly on a $G$-neighborhood of the image. In \cite{e2} the linearity assumption on the group $G$ is dropped and the assumptions on the space $X$ weakened, see Theorem~\ref{erik}. 

Here we will prove the following:\\

{\bf Theorem~\ref{emb}} \emph{Let $G$ be a Lie group and let $X \in \M$. Then there exists a $G$-embedding $e\colon X \rightarrow L$ where $L$ is a locally convex linear $G$-space such that $e(X)$ is a closed subset of some $G$-invariant convex subset $C$ of $L$ and $e(X)$ has some $G$-neighborhood $V$ in $C$ such that $V \in \M$.}\\

Using Theorem~\ref{emb} we prove that any $G$-ANE-$\M$ has the $G$-homotopy type of some $G$-CW complex, and that the classes $G$-ANE and $G$-ANR are the same.

\section{Preliminaries}

Throughout the paper $G$ will denote an arbitrary Lie group unless otherwise is stated, where Lie groups are defined to be Hausdorff and second countable. All $G$-spaces are assumed Hausdorff.

A completely regular $G$-space $X$ is said to be \emph{Cartan} if every point has a neighborhood $V$ such that the closure of the set $\{g \in G:gV \cap V \neq \emptyset\}$ is compact.

The action of $G$ on a completely regular $G$-space $X$ is \emph{proper} if for any pair of points $x,y \in X$ there exist neighborhoods $V_x$, $V_y$ of $x$ and $y$ such that the closure of the set $\{g \in G\colon gV_x \cap V_y \neq \emptyset\}$ is compact. We say that $X$ is a \emph{proper $G$-space}. The action of $G$ on a completely regular $G$-space $X$ is \emph{Palais proper} if for any point $x \in X$ there exists a neighborhood $V_x$ such that any point $y \in X$ has a neighborhood $V_y$ for which the closure of $\{g \in G:gV_x \cap V_y \neq \emptyset\}$ is compact. Then we say that $X$ is a \emph{Palais proper $G$-space}.

Clearly a Palais proper $G$-space is proper, and any proper $G$-space must be a Cartan $G$-space.

Let $H$ be a closed subgroup of $G$. A subset $S$ of a $G$-space $X$ is an \emph{$H$-slice} if $GS$ is open in $X$ and there exists a $G$-map $f \colon GS \rightarrow G/H$ such that $S=f^{-1}(eH)$. The set $S$ is a \emph{slice at the point $x \in X$} if $x \in S$ and $S$ is a $G_x$-slice.

By \cite[Theorem~2.3.3]{pa2} we know that in a Cartan $G$-space, there exists a slice at every point and the isotropy subgroup of $G$ at any point is compact.

A $G_x$-slice is also a $G_{gx}$-slice by the following lemma:

\begin{lem} \label{conj}
Suppose that $X$ is a Cartan $G$-space and let $S$ be a slice at $x \in X$. Then $S$ is an $H$-slice for any subgroup $H$ of $G$ which is conjugate to $G_x$.
\end{lem}

\begin{proof}
Suppose that $H=gG_xg^{-1}$ for some $g \in G$. By \cite[Proposition~1.1.5]{pa2} we have a map $h\colon G/G_x \approx_G Gx=Ggx \approx_G G/G_{gx}=G/gG_xg^{-1}=G/H$. Hence $h \circ f \colon GS \rightarrow G/G_x \circ G/H$ is the map associated with $S$ as an $H$-slice.
\end{proof}

The following characterization of slices is extremely useful.

\begin{thm} {\rm \cite[Theorem 2.1.4]{pa2}} \label{slicechar}
Suppose $X$ is a Cartan $G$-space and let $H$ be a compact subgroup of $G$. A subset $S \subset X$ is an $H$-slice if and only if the following conditions hold:
\begin{itemize}
\item[i)] $S$ is closed in $GS$,
\item[ii)] $S=HS$,
\item[iii)] $gS \cap S \neq \emptyset$ implies $g \in H$,
\item[iv)] $GS$ is open in $X$,
\item[v)] $S$ has a neighbohood $V$ in $GS$ such that the closure of $\{g \in G:gV \cap V \neq \emptyset\}$ is compact.
\end{itemize}
\end{thm}

We say that the open set $GS$ is a \emph{tubular neighborhood} (of $x$) if $S$ is an $H$-slice (a slice at $x$).

\begin{defn}[Tubular covering]
A \emph{tubular covering} of a $G$-space $X$ is a covering of $X$ by tubular neighborhoods.
\end{defn}

\begin{lem} \label{refine}
An open $G$-subset of a tubular neighborhood is a tubular neighborhood. Hence an open refinement by $G$-sets of a tubular covering is a tubular covering.
\end{lem}

\begin{proof}
Suppose $H$ is a closed subgroup of $G$, let $S$ be an $H$-slice and let $f\colon GS \rightarrow G/H$ be a corresponding $G$-map. Let $U$ be an open $G$-subset of $GS$. Now $S'=S \cap U$ is an $H$-slice since $GS'=U$ is open and $f'=f|U \rightarrow G/H$ is a $G$-map where $(f')^{-1}(eH)=f^{-1}(eH) \cap U=S \cap U=S'$.
\end{proof}

\begin{lem} \label{union}
If $\{GS_k\}_{k \in K}$ is a family of disjoint tubular neighborhoods in $X$, where $S_k$ is an $H$-slice for all $k \in K$, then $S=\bigcup_{k \in K} S_k$ is an $H$-slice.
\end{lem}

\begin{proof}
Here a corresponding map $f \colon GS=\bigcup_{k \in K} GS_k \rightarrow G/H$ is defined by $f|GS_k=f_k$ where $f_k\colon GS_k \rightarrow G/H$ is a map corresponding to the slice $S_k$.
\end{proof}

Given a closed subgroup $H$ of $G$ and an $H$-space $S$ there is an action of $H$ on the product $G \times S$ given by $h(g,s)=(gh^{-1}, hs)$. We denote by $G \times_H S$ the quotient space $(G \times S)/H$, which is called the \emph{twisted product} of $G$ and $S$ with respect to $H$.

\begin{prop} \cite[Proposition 1.18]{e1}
Let $H$ be a closed subgroup of $G$ and let $S$ be an $H$-slice in a $G$-space $X$. Then
$$
G \times_H S \approx_G GS.
$$
\end{prop} 

We will denote by $\M$ the class of proper metrizable $G$-spaces $X$ which have a metrizable orbit space $X/G$. By \cite[Theorem B]{anne}, for a proper metrizable $G$-space $X$, the metrizability of $X/G$ is equivalent to $X/G$ being paracompact, or $X$ admitting a $G$-invariant metric.

\section{Countability of tubular coverings}

The following lemma makes use of a technique originating with J.~Milnor, see \cite[Theorem 1.8.2]{pa1}.

\begin{lem} \label{thm}
Let $G$ be a Lie group, let $X$ be a Cartan $G$-space and suppose that $X/G$ is paracompact\footnote{We assume that paracompact spaces are Hausdorff.}. Then $X$ admits a countable and locally finite tubular covering.
\end{lem}

\begin{proof}
By \cite[Proposition 3.6]{an} we know that $G$ has at most countably many compact conjugacy types, represented by compact subgroups $H_n$ of $G$ where $n \in \N$. Suppose that $\{GS_i\}_{i \in I}$ is a tubular covering of $X$ such that for each $i \in I$, $S_i$ is a slice at some point $x_i \in X$, where $G_{x_i}$ is a compact subgroup of $G$. By Lemma~\ref{conj} we may assume that each $S_i$ is a $H_n$-slice for some $n \in \N$. Let $\{\varphi_i\colon X \rightarrow [0,1]\}_{i \in I}$ be a $G$-invariant partition of unity subordinate to $\{GS_i\}_{i \in I}$. Such a partition of unity exists because $X/G$ is paracompact.

For each finite $T \subset I$, denote ´
$$
W(T)=\{x \in X: \varphi_i(x)>\varphi_j(x) \textrm{ for all } i \in T \textrm{ and for all } j \in I \setminus T\}.
$$
Denote by $u_T \colon X \rightarrow [0,1]$ the continuous $G$-invariant map
$$
u_T(x)=\max \{0, \min_{i \in T, j \in I \setminus T} (\varphi_i(x)-\varphi_j(x))\}.
$$
Then $W(T)=u_T^{-1}(0,1]$ is open and $G$-invariant in $X$.

If $x \in W(T)$ then $\varphi_i(x) > \varphi_j(x) \ge 0$ for all $i \in T, j \in I \setminus T$ so in particular $x \in \varphi_i^{-1}(0,1]$ for all $i \in T$. Thus $W(T)$ is an open $G$-invariant subset of $\varphi_i^{-1}(0,1] \subset GS_i$ for each $i \in T$, hence $W(T)$ is a tubular neighborhood by Lemma~\ref{refine}.

If $\textrm{Card }T=\textrm{Card }T'$ and $T \neq T'$ then $W(T) \cap W(T')=\emptyset$, because if $i \in T \setminus T'$, $j \in T' \setminus T$ and $x \in W(T) \cap W(T')$ then simultaneously $\varphi_i(x)>\varphi_j(x)$ and $\varphi_j(x)<\varphi_i(x)$, which is impossible.

Define 
$$
W^m_n=\bigcup \{W(T):\textrm{Card }T=n \textrm{ and }W(T) \textrm{ is an } H_m \textrm{-tubular neighborhood}\}
$$
for all $n,m \in \N$. Now $W^m_n$ is an $H_m$-tubular neighborhood by Lemma~\ref{union} because it is a disjoint union of $H_m$-tubular neighborhoods, and $\{W^m_n\}_{(m,n) \in \N \times \N}$ is a countable tubular covering of $X$. Denote by $\pi_X \colon X \rightarrow X/G$ the projection onto the orbit space; since $X/G$ is paracompact, the open covering $\{\pi_X(W^m_n)\}_{(m,n) \in \N \times \N}$ of $X/G$ admits a precise locally finite refinement by \cite[VIII 1.4]{dug}; in particular this refinement is countable. Denote it $\{V_n\}_{n \in \N}$; now $\{\pi_X^{-1}(V_n)\}_{n \in \N}$ is a countable locally finite refinement of $\{W^m_n\}_{(m,n)\in \N \times \N}$ by $G$-sets; hence it is a countable locally finite tubular covering of $X$ by Lemma~\ref{refine}, and we are done.
\end{proof}

\section{Homeomorphism properties of isovariant $G$-maps}
 
The main result in this section is an important homeomorphism property of isovariant $G$-maps between Cartan $G$-spaces, which we will use in proving our main theorem. Recall that a map $f\colon \to Y$ between $G$-spaces is isovariant if $G_x=G_{f(x)}$ for all $x \in X$. Our lemma builds on the following result for compact transformation groups:

\begin{lem} \label{open} {\rm \cite[Exercise~10 of Chapter~I]{bre}}
Let $H$ be a compact group,\footnote{All topological groups are assumed to be Hausdorff.} and let $f \colon X \rightarrow Y$ be an isovariant $H$-map between $H$-spaces. Now $f$ is an open map if and only if the induced map $\bar{f}\colon X/H \rightarrow Y/H$ is an open map.
\end{lem}

\begin{lem} \label{open2}
Suppose that $X$ and $Y$ are Cartan $G$-spaces and that $f\colon X \rightarrow Y$ is an isovariant $G$-map. Then $f$ is a $G$-homeomorphism if and only if the induced map $\bar{f}\colon X/G \rightarrow Y/G$ is a homeomorphism.
\end{lem}

\begin{proof}
We have a commutative diagram

$$
\begin{diagram}
\node{X} \arrow{e,t}{f} \arrow{s,l}{\pi_X} \node{Y} \arrow{s,r}{\pi_Y}\\
\node{X/G} \arrow{e,b}{\bar{f}} \node{Y/G}
\end{diagram}
$$

and since $\pi_X$ is continuous and $\pi_Y$ is open, it is clear that if $f$ is open then so is $\bar{f}$. Hence if $f$ is a homeomorphism, so is $\bar{f}$.

Now suppose $\bar{f}$ is a homeomorphism; it then easily follows from the bijectivity of $\bar{f}$ and the fact that $f$ restricts to a bijection between orbits,\footnote{Any isovariant $G$-map restricts to a bijection on the orbits.} that $f$ is a bijection, and it suffices to show that $f$ is open. Being Cartan, $Y$ has a covering $\{GS^Y_i\}_{i \in I}$ where each $S^Y_i$ is a slice at a point $y_i \in Y$ such that the isotropy subgroups $G_{y_i}$ are compact for all $i \in I$. Let $i \in I$. Then $S^X_i=f^{-1}(S^Y_i)$ is a slice at $x_i=f^{-1}(y_i)$. Note that by isovariance $G_{y_i}=G_{x_i}$. We consider now only the restrictions
$$
\begin{diagram}
\node{GS^X_i} \arrow{e,t}{f} \arrow{s,l}{\pi_X} \node{Gf(S^X_i)=GS^Y_i} \arrow{s,r}{\pi_Y}\\
\node{GS^X_i/G} \arrow{e,b}{\bar{f}} \node{Gf(S^X_i)/G=GS^Y_i/G}
\end{diagram}
$$
Since $GS^X_i/G \approx S^X_i/G_{x_i}$ and $GS^Y_i/G \approx S^Y_i/G_{x_i}$ by \cite[Lemma 1.23]{e1}, we see that the restriction $f|S^X_i\colon S^X_i \rightarrow S^Y_i$ induces a map $\tilde{f}\colon S^X_i/G_{x_i} \rightarrow S^Y_i/G_{x_i}$ given by $\tilde{f}(G_{x_i}s)=G_{x_i}(f(s))$, which makes the following diagram commute:
$$
\begin{diagram}
\node{S^X_i} \arrow{e,t}{f} \arrow{s,l}{\pi} \node{S^Y_i} \arrow{s,r}{\pi'}\\
\node{S^X_i/G_{x_i}} \arrow{e,b}{\tilde{f}} \arrow{s,l}{\approx} \node{S^Y_i/G_{x_i}} \arrow{s,r}{\approx}\\
\node{GS^X_i/G} \arrow{e,b}{\bar{f}} \node{GS^Y_i/G}
\end{diagram}
$$
Hence since $\bar{f}$ is open, so is $\tilde{f}$ and thus by Lemma~\ref{open} the restriction $f|S^X_i$ is open since $G_{x_i}$ is compact by assumption.

The map $\id_G \times f|S^X_i$ induces a map $G \times_{G_{x_i}} f|S^X_i\colon G \times_{G_{x_i}} S^X_i \rightarrow G \times_{G_{x_i}} S^Y_i$ as in \cite[Proposition II 2.1]{bre} which makes the following diagram commute:
$$
\begin{diagram}
\node{G \times S^X_i} \arrow{e,t}{\id_G \times f|S^X_i} \arrow{s} \node{G \times S^Y_i} \arrow{s}\\
\node{G \times_{G_{x_i}} S^X_i} \arrow{e,t}{G \times_{G_{x_i}} f|S^X_i} \arrow{s,l,r}{[g,s] \mapsto gs, \ \approx_G} \node{G \times_{G_{x_i}} S^Y_i} \arrow{s,r}{[g,y] \mapsto gy, \ \approx_G}\\
\node{GS^X_i} \arrow{e,b}{f} \node{GS^Y_i}
\end{diagram}
$$
and since $f|S^X_i$ is open, $G \times_{G_{x_i}} f|S^X_i$ is open by \cite[Proposition II 2.1]{bre}. Hence $f|GS^X_i$ is open onto its image, i.e.~it is a homeomorphism onto its image.

Since the $S^Y_i$ are slices, the set $GS^Y_i$ is open in $Y$ for each $i \in I$. Hence $f$ is open, i.~e. it is a homeomorphism.
\end{proof}

\section{Embedding theorem}

Now we are able to prove the following theorem, which is the main result of this note:

\begin{thm} \label{emb}
Let $G$ be a Lie group and let $X \in \M$. Then there exists a $G$-embedding $e\colon X \rightarrow L$ where $L$ is a locally convex linear $G$-space such that $e(X)$ is a closed subset of some $G$-invariant convex subset $C$ of $L$ and $e(X)$ has some $G$-neighborhood $V$ in $C$ such that $V \in \M$.
\end{thm}

\begin{defn}[LCL $G$-space]
A locally convex linear $G$-space (for short, an LCL $G$-space) is a $G$-space $L$ which is a locally convex topological vector space where each element $g \in G$ represents a linear map $L \rightarrow L$.
\end{defn}

We will prove Theorem~\ref{emb} by finding an isovariant $G$-map $e \colon X \rightarrow L$ which induces an embedding $X/G \rightarrow L/G$ in suitable conditions, and by applying Lemma~\ref{open2}. The same method has been used by J.~Jaworowski \cite{jaw82} and of G.~Bredon \cite[Chapter II.10]{bre} for compact Lie group actions.

\begin{lem} \label{prodlcl}
A product of LCL $G$-spaces with diagonal action is an LCL $G$-space.
\end{lem}

\begin{proof}
Let $X_i$ be a family of LCL $G$-spaces where the indices $i$ run through some set $I$. Set $X=\prod_{i \in I} X_i$; now it is clear that $X$ is a topological vector space and a $G$-space where each $g \in G$ represents a linear map $X \rightarrow X$. Furthermore $X$ is locally convex:\\

Suppose $U$ is an open neighborhood of $x \in X$. Then there is some basis element $V=\prod_{i \in I} V_i$ such that $x \in V \subset U$. For $V$ to be an element in the basis, we know that each $V_i$ is open in $X_i$ and $V_i=X_i$ for all but finitely many $i \in I$. Now for each $V_i \neq X_i$ there exists a convex neighborhood $W_i$ of $\pr_i(x)$ such that $W_i \subset V_i$. Whenever $V_i=X_i$, set $W_i=X_i$. Define $W=\prod_{i \in I} W_i$. Now $W$ is a convex neighborhood of $x$ and $W \subset V \subset U$. 
\end{proof}

We are going to use the following result, which is obtained in \cite[Section~3]{e2}, although there it is not explicitly stated as a theorem:

\begin{thm} \label{erik}
Let $G$ be a Lie group and assume that a $G$-space $X \in \M$ admits a metric $d$ which satisfies
$$
(\ast) \quad \forall \ r>0 \ \forall \ x \in X\colon \bar{B}^d_r(x) \textrm{ is compact.} 
$$

Then there exists a Banach $G$-space $B$ and a closed $G$-embedding $i \colon X \rightarrow B$ such that $i(X) \subset B \setminus \{\bar{0}\}$ and $G$ acts properly on $B \setminus \{\bar{0}\}$.
\end{thm}

\begin{rem}
\begin{itemize}
\item[i)] A normed vector space is locally convex.
\item[ii)] From the construction of the space $B$ one easily sees that the norm in $B$ is $G$-invariant, inducing a $G$-invariant metric on $B$. That is, $B \setminus \{\bar{0}\} \in \M$.
\end{itemize}
\end{rem}

Suppose that the $G$-space $X$ has a global $H$-slice $S$, i.e.~$X=GS$. Then there exists an isovariant $H$-map 
$$
\varphi\colon S \rightarrow \prod_{n \in \N} E_n
$$
where each $E_n$ is a Euclidean representation space for $H$ (see \cite[Lemma 5]{an}). The twisted product $G \times_H E_n$ is a $G$-space and there is an $H$-embedding $i_n\colon E_n \rightarrow G \times_H E_n$ defined by $i_n(x)=[e,x]$. This defines an isovariant $H$-map
$$
\tilde{\varphi}\colon S \stackrel{\varphi}{\rightarrow} \prod_{n \in \N} E_n \xrightarrow{\prod i_n} \prod_{n \in \N} G \times_H E_n.
$$
Using this we obtain a $G$-map $\psi\colon X=GS \rightarrow \prod_{n \in \N} G \times_H E_n$ by setting $\psi(gx)=g\tilde{\varphi}(x)$.

\begin{lem}
The map $\psi$ is isovariant.
\end{lem}

\begin{proof}
By equivariance $G_{gs} \subset G_{\psi(gs)}$ for all $gs \in GS$. Thus we only need to show ''$\supset$''.\\

Assume that $g\tilde{\varphi}(S) \cap \tilde{\varphi}(S) \neq \emptyset$ for some $g \in G$. Then $[g, \varphi_n(s_1)]=[e,\varphi_n(s_2)]$ for all $n \in \N$ and for some $s_1,s_2 \in S$, where $\varphi_n=\pr_n \circ \varphi$. But then $g \in H$.

If we now assume that $\bar{g} \in G_{\psi(gs)}$ for some $gs \in GS$, then $\bar{g}\psi(gs)=\psi(gs)$, giving $\bar{g}g\tilde{\varphi}(s)=g\tilde{\varphi}(s)$, hence $g^{-1}\bar{g}g\tilde{\varphi}(s)=\tilde{\varphi}(s)$ and thus by the previous argument $g^{-1}\bar{g}g \in H_{\tilde{\varphi}(s)}=H_s$.

But then $g^{-1}\bar{g}gs=s$, hence $\bar{g}gs=gs$ and thus $\bar{g} \in G_{gs}$. It follows that $G_{\psi(gs)}=G_{gs}$ for all $gs \in GS$ and thus $\psi$ is isovariant.
\end{proof}

Now $G \times_H E_n$ is a proper and $N_2$ manifold (manifold because $G \times_H E_n \rightarrow G/H$ is a vector bundle by \cite[Theorem~2.26]{kaw} and proper by \cite[Proposition~1.3]{e1}). It has a metric with the property~$(\ast)$ of Theorem~\ref{erik} because, being separable metric and finite-dimensional, $G \times_H E_n$ can be embedded in some Euclidean space. Furthermore, $(G \times_H E_n)/G$ is metrizable by \cite[Theorem~4.3.4]{pa2}, giving $G \times_H E_n \in \M$. Hence, by Theorem~\ref{erik} we obtain a $G$-embedding $G \times_H E_n \rightarrow B_n$, where $B_n$ is a Banach $G$-space and $G$ acts properly on $B_n \setminus \{\bar{0}\}$. This gives an isovariant $G$-map
$$
\tilde{\psi}\colon X=GS \stackrel{\psi}{\rightarrow} \prod_{n \in \N} G \times_H E_n \rightarrow \prod_{n \in \N} B_n =: \tilde{Z},
$$
where $\tilde{\psi}(GS) \subset \tilde{Z} \setminus \{\bar{0}\}$.

\begin{lem} \label{0.10}
$G$ acts properly on $\tilde{Z} \setminus \{\bar{0}\}$.
\end{lem} 

\begin{proof}
$G$ acts properly on $B_n \setminus \{\bar{0}\}$ for each $n \in \N$ and hence the space $\tilde{Z}_n=\left(\prod_{i=1}^{n-1} B_i \right) \times (B_n \setminus \{\bar{0}\}) \times \left(\prod_{i=n+1}^\infty B_i \right)$ is a proper $G$-space for all $n \in \N$. Furthermore $\tilde{Z} \setminus \{\bar{0}\} = \bigcup_{n \in \N} \tilde{Z}_n$. Now $\tilde{Z}$ is a Cartan $G$-space since any completely regular $G$-space which is the union of open Cartan $G$-subspaces is Cartan. Furthermore, there is a $G$-invariant metric $d$ on $\tilde{Z}$, which induces a pseudometric $\bar{d}$ on $\tilde{Z}/G$ by $\bar{d}(\tilde{x},\tilde{y})=d(Gx,Gy)$, where $\tilde{x}=\pi(x)$ and $\pi\colon \tilde{Z} \rightarrow \tilde{Z}/G$. The pseudometric $\bar{d}$ induces the quotient topology on $\tilde{Z}/G$ since 
$$
y \in B_d(x,r) \Rightarrow \bar{d}(\tilde{x},\tilde{y}) \le d(x,y)<r \Rightarrow \pi(y) = \tilde{y} \in B_{\bar{d}}(\tilde{x},r)
$$
and
$$
\begin{array}{ll}
\tilde{y} \in B_{\bar{d}}(\tilde{x},r) & \Rightarrow \inf_{g \in G} d(x,gy)=\bar{d}(\tilde{x},\tilde{y})<r \Rightarrow gy \in B_d(x,r) \textrm{ for some } g \in G\\
 & \Rightarrow \tilde{y}=\pi(gy) \in \pi B_d(x,r)
\end{array}
$$
imply $\pi B_d(x,r)=B_{\bar{d}}(\tilde{x},r)$.

Hence $(\tilde{Z} \setminus \{\bar{0}\})/G$ is pseudometrizable, and thus regular. It follows that the action of $G$ on $\tilde{Z} \setminus \{\bar{0}\}$ is Palais proper by \cite[Proposition 1.2.5]{pa2}, so in particular it is proper.
\end{proof}

Next we pass from the global slice situation above to the more complicated situation with arbitrary spaces from $\M$.

\begin{lem} \label{4.8}
Let $X$ be as in Theorem~\ref{emb}. There exists an LCL $G$-space $Z$ and an isovariant $G$-map $\tilde{f}\colon X \rightarrow Z$ such that $\tilde{f}(X) \subset Z \setminus \{\bar{0}\}$ and $G$ acts properly on $Z \setminus \{\bar{0}\}$.
\end{lem}

\begin{proof}
By Lemma~\ref{thm} we may assume that $\{GS_n\}_{n \in \N}$ is a locally finite tubular covering of $X$. By Lemma~\ref{refine} and by the normality of $X/G$ we may let $\{GR_n\}_{n \in \N}$ and $\{GT_n\}_{n \in \N}$ be similar coverings such that $\overline{GT}_n \subset GR_n \subset \overline{GR}_n \subset GS_n$ for each $n \in \N$. According to the discussion above each $GS_n$ admits an isovariant $G$-map $\tilde{f}_n$ into an LCL $G$-space $Z_n$, where $\tilde{f}_n(GS_n) \subset Z_n \setminus \{\bar{0}\}$, and $G$ acts properly on $Z_n \setminus \{\bar{0}\}$. For each $n \in \N$ we may construct a $G$-invariant map $\lambda_n\colon X \rightarrow I$ such that $\lambda_n(\overline{GT}_n)=\{1\}$ and $\lambda_n(X \setminus GR_n)=\{0\}$, since $X/G$ is normal.

We obtain a continuous and isovariant $G$-map $\tilde{f}\colon X \rightarrow \prod_{n \in \N} Z_n=:Z$ by setting $\tilde{f}(x)=(\lambda_1(x)\tilde{f}_1(x),\ldots, \lambda_i(x)\tilde{f}_i(x),\ldots)$.

Now $Z$, being a product of LCL $G$-spaces, is an LCL $G$-space by Lemma~\ref{prodlcl}.

Clearly $\tilde{f}(X) \subset Z \setminus \{\bar{0}\}$, and the proof that $G$ acts properly on $Z \setminus \{\bar{0}\}$ goes just as in Lemma~\ref{0.10}.
\end{proof}

\begin{lem}
With $X$ as in Theorem~\ref{emb}, there exists an embedding $e\colon X \rightarrow L$ where $L$ is an LCL $G$-space, $e(X) \subset C$ for some convex $G$-subset $C$ of $L$ and $e(X)$ has a $G$-neighborhood $V$ in $C$ such that $V \in \M$.
\end{lem}

\begin{proof}
Let $h\colon X/G \rightarrow B$ be an embedding of the metrizable space $X/G$ into a Banach space $B$ such that $h(X/G)$ is a closed subset of $C'$ where $C'$ is a convex subset of $B$. Such a map exists by the Wojdys{\l}awski embedding theorem. With $Z$ as in Lemma~\ref{4.8}, define a map $e\colon X \rightarrow Z \times B=L$ by setting $e=(\tilde{f},h \circ \pi_X)$ where $\pi_X\colon X \rightarrow X/G$ is the natural projection, and let $G$ act trivially on $B$. Now we have

\begin{itemize}
\item[i)] $e$ is an isovariant $G$-map because $\tilde{f}$ is isovariant.
\item[ii)] $e$ is injective since the induced map $\bar{e}\colon X/G \rightarrow Z/G \times B$ is injective and $e$ is isovariant. The map $\bar{e}$ is injective because $h$ is injective.
\item[iii)] The map $\bar{e}$ is a homeomorphism onto its image, whose inverse $\bar{e}^{-1}\colon \bar{e}(X/G)=\pi_Z(\tilde{f}(X)) \times h(X/G) \rightarrow X/G$ is given by $h^{-1} \circ \pr_2$.
\item[iv)] $e(X) \subset Z \times C'=C$ which is convex, and $G$ acts properly on $Z\setminus \{\bar{0}\} \times C'=V$, which is a neighborhood of $e(X)$ in $C$.
\end{itemize}

Now $e(X) \subset V$ which is proper, and in particular $e(X)$ is Cartan. Thus by Lemma~\ref{open2} the map $e$ is a homeomorphism onto its image. The space $V/G=(Z \setminus \{\bar{0}\} \times C')/G$ is metrizable since $Z \setminus \{\bar{0}\} \times C'$ is proper and admits a $G$-invariant metric. Thus $V \in \M$.
\end{proof}

\begin{lem}
The image $e(X)$ is a closed subset of $C=Z \times C'$.
\end{lem}

\begin{proof}
Assume that $y=(u,v) \in (Z \times C') \setminus e(X)$. In case $v \notin h(X/G)$ then since $h$ is a closed embedding, there exists a neighborhood $U$ of $v$ in $C'$ such that $U \cap h(X/G)=\emptyset$. Thus $Z \times U$ is a neighborhood of $(u,v)$ which does not intersect $e(X)$.

In case $v \in h(X/G)$ we have $v=h(\pi_X(x))$ for some $x \in X$ and then we must have $Gu \cap G\tilde{f}(x)= \emptyset$. (If not then $u=g\tilde{f}(x)=\tilde{f}(gx)$ for some $g \in G$, giving $v=h(\pi_X(x))=h(\pi_X(gx))$ which implies $y=(u,v)=e(gx) \in e(X)$.)

If $u \neq \bar{0}$ then $u$ and $\tilde{f}(x)$ are points in the open subset $Z\setminus \{\bar{0}\}$ of $Z$, where $(Z\setminus \{\bar{0}\})/G$ is metrizable. Thus there exist disjoint $G$-neighborhoods $U$ and $W$ of $Gu$ and $G\tilde{f}(x)$ in $Z\setminus\{\bar{0}\}$ and hence in $Z$. If $u=\bar{0}$ then $\tilde{u}$ is closed in $Z/G$ which is pseudometrizable and hence regular, so there exist again disjoint $G$-neighborhoods $U$ and $W$ of $Gu$ and $G\tilde{f}(x)$ in $Z$. Thus we see that in any case there exist disjoint $G$-neighborhoods $U$ and $W$ of $Gu$ and $G\tilde{f}(x)$ in $Z$. 

Since $\tilde{f}$ is continuous and equivariant there exists a $G$-neighborhood $\tilde{W}$ of $Gx$ in $X$ such that $\tilde{f}\tilde{W} \subset W$. Since $\pi_X$ is an open map and $h$ is an embedding there exists an open neighborhood $M$ of $h(\pi_X(x))$ in $C'$ such that $M \cap h(\pi_X(X))=h(\pi_X(\tilde{W}))$.

Clearly $U \times M$ is an open neighborhood of $y$ in $Z \times C'$; we show that it is disjoint from $e(X)$. If $(x,w) \in e(X) \cap U \times M$ then $x=\tilde{f}(z) \in U$ and $w=h(\pi_X(z)) \in M$ for some $z \in X$. But then $\tilde{f}(z) \notin W$ giving $z \notin \tilde{W}$ and since $\tilde{W}$ is $G$-invariant we have $\pi_X(z) \notin \pi_X(\tilde{W})$ which gives $h(\pi_X(z)) \notin h(\pi_X(\tilde{W}))$ since $h$ is injective.

However, $z \in X$ and $h(\pi_X(z)) \in M$ implies that $h(\pi_X(z)) \in M \cap h(\pi_X(X))=h(\pi_X(\tilde{W}))$ which gives a contradiction. Hence we must have $e(X) \cap (U \times M)=\emptyset$.

It follows that $e(X)$ is closed in $Z \times C'$.
\end{proof}

\begin{proof}[The rest of the proof of Theorem~\ref{emb}] If we now set $L=Z \times B$, then $L$ is an LCL $G$-space. Set $C=Z \times C'$ and $V=Z \setminus \{\bar{0}\} \times C'$ as before. Then $V$ is a $G$-neighborhood of $e(X)$ in $C$, $C$ is a $G$-invariant convex subset of $L$, $e(X)$ is closed in $V$, and $V \in \M$. Hence the theorem is true. 
\end{proof}

\section{Applications}

\subsection{Any $G$-ANE is $G$-homotopy equivalent to some $G$-CW complex}

We say that a $G$-space $Y \in \M$ is a $G$-equivariant absolute neighborhood retract ($G$-ANE) for $\M$, written $Y \in G$-ANE-$\M$, if for any $G$-space $X \in \M$ and any closed $G$-invariant subset $A$ in $X$ with a $G$-map $f \colon A \rightarrow Y$ there exists a $G$-extension $F \colon U \rightarrow Y$ of $f$ over some $G$-neighborhood $U$ of $A$ in $X$.

It has been shown by M.~Murayama \cite{mu}\footnote{There is an error in the equivariant embedding theorem used in \cite{mu}, since the action on the linear space is not necessarily continuous. See \cite{an*} for a counterexample and a correct equivariant embedding theorem.} that for a compact Lie group $G$, any $X \in G$-ANE-$\M$ has the $G$-homotopy type of some $G$-CW complex. In \cite{e1} E.~Elfving showed that, given a linear Lie group $G$, any locally linear manifold $M \in \M$ (then $M \in G$-ANE-$\M$ as well) has the $G$-homotopy type of some $G$-CW complex, and in \cite{e2} he removes the linearity assumption on the Lie group $G$. Using Theorem~\ref{emb} we get

\begin{cor}
Any $G$-ANE for $\M$ is $G$-homotopy equivalent to some $G$-CW complex.
\end{cor}

\begin{proof}
Suppose that $X \in G$-ANE-$\M$. Then $X \in \M$ by definition, and according to the previous embedding theorem we may assume that $X$ is a closed $G$-subset of a convex $G$-subset $C$ of an LCL $G$-space $L$, where $X$ has some $G$-neighborhood $V$ in $C$ such that $V \in \M$. Then there exists a $G$-extension $r\colon U \rightarrow X$ of $\id\colon X \rightarrow X$, where $U$ is a $G$-neighborhood of $X$ in $V$, and hence in $C$. Here $C$ is locally convex, being a convex subset of a locally convex space, and $r$ is a $G$-retraction. Now one can use the argument of \cite[Proposition 11, Theorem 5 and Theorem 6]{e2} to show that $X$ has the $G$-homotopy type of some $G$-CW complex.
\end{proof}

\subsection{Equivalence of the classes $G$-ANE-$\M$ and $G$-ANR-$\M$}

We say that a $G$-space $X \in \M$ is a $G$-equivariant absolute neighborhood retract ($G$-ANR) for $\M$, written $X \in G$-ANR-$\M$, if, whenever there exists a closed $G$-embedding $i\colon X \rightarrow Y$ of $X$ into some $G$-space $Y \in \M$, then there exists a $G$-neighborhood retraction $r \colon U \rightarrow i(X)$ where $U$ is a $G$-neighborhood of $i(X)$ in $Y$.

It is easy to show that that $G$-ANE-$\M \subset G$-ANR-$\M$. Here we show that $G$-ANR-$\M \subset G$-ANE-$\M$, following the classical proof from the non-equivariant case (see, for instance, \cite{hu}). 

Suppose that $Y \in G$-ANR-$\M$. Now by Theorem~\ref{emb} we may assume that $Y$ is a closed $G$-subset of a convex $G$-subset $C$ of some LCL $G$-space, where $Y$ has a $G$-neighborhood $U$ in $C$ such that $U \in \M$. Since $Y \in G$-ANR-$\M$, there exists a $G$-neighborhood $V$ of $Y$ in $U$ (and hence in $C$) and a $G$-retraction $r\colon V \rightarrow Y$. Let $X \in \M$ and let $A$ be a closed $G$-invariant subset of $X$ with a $G$-map $f \colon A \rightarrow Y$. By the equivariant Dugundji extension theorem \cite[Corollary 1]{an2} we know that $C \in G$-ANE-$\M$; hence the map $i \circ f \colon A \rightarrow Y \hookrightarrow C$ admits a $G$-extension $h \colon W \rightarrow C$, where $W$ is a $G$-neighborhood of $A$ in $X$. Set $W'=W \cap h^{-1}(V)$. Now $W'$ is a $G$-neighborhood of $A$ in $X$ and the map $r \circ h|W' \colon W' \rightarrow V \rightarrow Y$ is a neighborhood $G$-extension of $f$. It follows that $Y \in G$-ANE-$\M$. 

\end{document}